\documentclass[a4paper,12pt]{amsart}
\usepackage{amsmath,amssymb}
\usepackage[left=2cm,right=2cm]{geometry}
\makeatletter
\@namedef{subjclassname@2020}{%
	\textup{2020} Mathematics Subject Classification}
\makeatother
\begin{document}

\title[A Schwarz-type lemma for squeezing function]{A Schwarz-type lemma for squeezing function on planar domains}

\author{Ahmed Yekta {\"O}kten}

\address{A. Y. {\"O}kten\\
	Institut de Math\'ematiques de Toulouse; UMR5219 \\
	Universit\'e de Toulouse; CNRS \\
	UPS, F-31062 Toulouse Cedex 9, France} \email{ahmed$\_$yekta.okten@math.univ-toulouse.fr}

	

\thanks{	The author received support from the University Research School EUR-MINT (State support managed
	by the National Research Agency for Future Investments program bearing the reference ANR-18-EURE-0023).}

\newcommand{\obj}[3]{\mathcal{F}^{#1}\mathbb{S}^{#2}\mathbf{G}_{#3}}

\newcommand{\dbar}{\bar\partial}
\newcommand{\genmat}{\lambda}
\newcommand{\polynorm}[1]{{|| #1 ||}}
\newcommand{\vnorm}[1]{\left\|  #1 \right\|}
\newcommand{\asspol}[1]{{\mathbf{#1}}}

\newcommand{\U}{{\mathscr{U}}}
\newcommand{\Uo}{\mathcal{U}(\Omega)}
\newcommand{\Fo}{\mathcal{F}(\Omega)}

\newcommand{\Uot}{\tilde{\mathcal{U}}(\Omega)}
\newcommand{\Cd}{\mathbb{C}^d}
\newcommand{\Cm}{\mathbb{C}^m}
\newcommand{\arO}{r_{\Omega,i,j}(z)}
\newcommand{\esO}{s_{\Omega,i,j}(\Omega)}
\newcommand{\ceOj}{C{z,j}(\Omega)}
\newcommand{\dist}{d}
\newcommand{\CN}{\mathbb{C}^N}
\newcommand{\CNp}{\mathbb{C}^{N^\prime}}
\newcommand{\Rd}{\mathbb{R}^d}
\newcommand{\Rn}{\mathbb{R}^n}
\newcommand{\RN}{\mathbb{R}^N}

\newcommand{\Ui}{\tilde{\mathbb{U}}_{\Gamma_i}(\Omega)}
\newcommand{\Ar}{\mathbb{A}_{r}}
\newcommand{\G}{\Gamma}
\newcommand{\N}{\mathbb{N}}
\newcommand{\dop}[1]{\frac{\partial}{\partial #1}}
\newcommand{\vardop}[3]{\frac{\partial^{|#3|} #1}{\partial {#2}^{#3}}}
\newcommand{\infnorm}[1]{{\left\| #1 \right\|}_{\infty}}
\newcommand{\onenorm}[1]{{\left\| #1 \right\|}_{1}}
\newcommand{\deltanorm}[1]{{\left\| #1 \right\|}_{\Delta}}
\newcommand{\omeganorm}[1]{{\left\| #1 \right\|}_{\Omega}}
\newcommand{\nequiv}{{\equiv \!\!\!\!\!\!  / \,\,}}
\newcommand{\bk}{\mathbf{K}}
\newcommand{\p}{\prime}
\newcommand{\tV}{\mathcal{V}}
\newcommand{\poly}{\mathcal{P}}
\newcommand{\ring}{\mathcal{A}}
\newcommand{\ringk}{\ring_k}
\newcommand{\ringktwo}{\mathcal{B}_\mu}
\newcommand{\germs}{\mathcal{O}}
\newcommand{\On}{\germs_n}
\newcommand{\mcl}{\mathcal{C}}
\newcommand{\formals}{\mathcal{F}}
\newcommand{\Fn}{\formals_n}
\newcommand{\dc}{\mathnormal{C}}
\newcommand{\autM}{{\Aut (M,0)}}
\newcommand{\autMp}{{\Aut (M,p)}}
\newcommand{\holmaps}{\mathcal{H}}
\newcommand{\biholmaps}{\mathcal{B}}
\newcommand{\autmaps}{\mathcal{A}(\CN,0)}
\newcommand{\jetsp}[2]{ G_{#1}^{#2} }
\newcommand{\njetsp}[2]{J_{#1}^{#2} }
\newcommand{\jetm}[2]{ j_{#1}^{#2} }
\newcommand{\glnc}{\mathsf{GL_n}(\C)}
\newcommand{\glmc}{\mathsf{GL_m}(\C)}
\newcommand{\glc}{\mathsf{GL_{(m+1)n}}(\C)}
\newcommand{\glk}{\mathsf{GL_{k}}(\C)}
\newcommand{\smC}{\mathcal{C}^{\infty}}
\newcommand{\anC}{\mathcal{C}^{\omega}}
\newcommand{\kC}{\mathcal{C}^{k}}
\def\dxdy{\delta_\Omega(x)\delta_\Omega(y)}
\def\B{\mathbb B}
\def\w{\omega_{\Omega}}
\def\C{\mathbb C}
\def\Cn{\mathbb{C}^n}
\def\el{L_e}
\def\kl{L_K}
\def\krl{L_{\kappa}}
\def\Bn{\mathbb{B}^n}
\def\R{\mathbb R}
\def\kr{\kappa_\Omega}
\def\kb{\hat{\kappa}}
\def\kbl{L_{\hat{\kappa}}}
\def\D{\mathbb D}
\def\L{\mathcal L}
\def\a{\alpha}
\def\endofproof{\hfill \square}
\def\d{\delta}
\def\g{g_\Omega}
\def\z{\zeta}
\def\ds{\displaystyle}
\def\sdd{\sqrt{\delta_{\Omega}(x)\delta_{\Omega}(y)}}
\def\jump{$\:$ \newline}

\newtheorem{Theorem}{Theorem}[section]
\newtheorem{cor}[Theorem]{Corollary}
\newtheorem{prop}[Theorem]{Proposition}
\newtheorem{lem}[Theorem]{Lemma}

\theoremstyle{definition}\newtheorem{Def}[Theorem]{Definition}

\theoremstyle{remark}
\newtheorem{Rem}[Theorem]{Remark}
\newtheorem{Exa}[Theorem]{Example}
\newtheorem{Exs}[Theorem]{Examples}

\numberwithin{equation}{section}



\begin{abstract}  
	With an easy application of maximum principle, we establish a Schwarz-type lemma
	for squeezing functions on finitely connected planar domains which bounds their squeezing functions with respect to invariants obtained from their canonical circularly slit representations. In particular, our result directly yields the explicit formulas
	of squeezing functions on doubly connected domains obtained by Ng, Tang and Tsai.
 \end{abstract}

\maketitle

	\section{Introduction}\label{int} 

Let $\Omega$ be a domain in $\Cn$ such that the set of holomorphic embeddings $\Uo$ of $\Omega$ into the unit ball $\Bn \subset \Cn$ is non-empty. The \emph{squeezing function} of $\Omega$, $S_{\Omega}:\Omega \longrightarrow (0,1] $ is defined as  \\
\[ S_{\Omega}(z):=\sup_{\substack{f\in\Uo\\f(z)=0}}\{r\in(0,1]:r\Bn \subset f(\Omega)\} \\
=\sup_{\substack{f\in\Uo\\f(z)=0}}\dist (0,\partial f(\Omega)). \]

The squeezing function is a biholomorphic invariant that arose from the study of invariant metrics on Teichmüller spaces of Riemann surfaces, see \cite{Y,LSY1,LSY2}. The papers \cite{YSK,DGZ1,DGZ2} provide an introduction to the topic. 

We follow the terminology of \cite{S}. By \textit{squeezing problem on $\Omega$} we mean the problem of determining $S_{\Omega}$ explicitly. A map $f\in\Uo$ is said to be \textit{extremal} for the squeezing problem for $z\in\Omega$ if $f(z)=0$ and $S_{\Omega}(z)=\dist(0,\partial f(\Omega))$. A normal families argument shows the existence of the extremal maps for the squeezing function problem. 

The boundary behaviour of the squeezing function is well-studied and exploited. See for instance \cite{DGZ1, DGZ2, F1, F2}. On the other hand, determining the squeezing function explicitly is a  hard problem, even in the case of planar domains. Recently, a non-trivial example is given in \cite{NTT} where the authors used a modified version of Löwner's differential equation in order to obtain the explicit formulas of squeezing functions of doubly connected planar domains. In \cite{GR} and \cite{S}, the same result is proven with different techniques and it is also established that biholomorphisms that take the domain onto circularly slit discs (see Section \ref{preliminaries} for the definition) in a suitable way are the only extremal maps for squeezing problems on doubly connected domains. 

In this paper, we prove Theorem \ref{mainresult}, which gives upper and lower bounds to the squeezing functions of finitely connected planar domains with respect to invariants obtained from their canonical circularly slit representations. Notably, our proof easily follows from an application of the maximum principle where we exploit the locally extremal properties of biholomorphisms that map a domain onto a circularly slit disc. Moreover, our method directly yields explicit formulas of squeezing functions of doubly connected domains and it gives the uniqueness result for extremal maps for this problem.

\section{Discussion and results}\label{preliminaries}

\subsection{Notation and statement of the main result} $\:$

Recall that a domain $\Omega\subset\C$ is said to be $n$-connected if $\mathbb{P}^1\setminus \Omega$ has $n$-components and it is said to have non-degenerate boundary if no component of $\mathbb{P}^1\setminus D$ is a point. Here $\mathbb P^1$ denotes the Riemann sphere. 

Let $\Delta$ denote the open unit disc in the complex plane. A domain $C\subset\Delta$ is said to be a circularly slit disc if the components of $\Delta \setminus C$ consists of closed proper arcs of circles centered at the origin. 

It is well-known that any multiply connected planar domain with non-degenerate boundary is biholomorphic to a circularly slit disc. Moreover, it is also known that any biholomorphism between two circularly slit discs fixing the origin and the unit circle must be a rotation. More explicitly, we have the following stronger result.

\begin{lem}\cite[Theorem 6.2]{C}
	Let $\Omega\subset\C$ be a finitely connected domain with non-degenerate boundary components $\{\G_k\}_{k \in \{1,2...,n\}}$ and $z\in \Omega$, for each ${j \in \{1,2...,n\}}$ there exists a unique biholomorphism $\psi_{z,j}$ taking $\Omega$ onto a circularly slit disc such that $\psi_{z,j}(z)=0$, $\psi'_{z,j}(z)=0$ and $\psi_{z,j}(\G_j)=\partial \Delta$ in the sense of boundary correspondence. 
\end{lem}

Note that boundary correspondence means that if $f:\Omega_1\rightarrow \Omega_2$ is a biholomophism, $\Gamma_1\subset\partial \Omega_1$ and $\Gamma_2\subset\partial \Omega_2$ we have $f(\Gamma_1)=\Gamma_2$ if for any sequence $\{z_k\}\subset \Omega_1$ we have $z_k$ tends to $\Gamma_1$ if and only if $f(z_k)\in\Omega_2$ tends to $\Gamma_2$. 

For a background on circularly slit discs, we refer the reader to \cite[Section 15.6]{C}.

Having this in mind, throughout the note we will use the following notation. As above let $\Omega\subset\mathbb C$ be an $n$-connected domain with non-degenerate boundary components $\{\G_k\}_{k \in \{1,2...n\}}$ and $z \in \Omega$. For $j \in \{1,2...n\}$ let $\psi_{z,j}$ denote the unique biholomorphism mapping $\Omega$ onto the circularly slit disc $C_{z,j}(\Omega)$ such that ${\psi_{z,j}}(z)=0$, ${\psi_{z,j}}(\G_j)=\partial \Delta$ and normalized so that $\psi_{z,j}'(z)>0$. For $i\neq j$ denote the components of $\Delta\setminus C_{z,j}(\Omega)$ as $\esO={\psi_{j,z}}(\G_i)$ and let $\arO=\dist(0,\esO)$. 

With this notation, our main result is the following. 
\begin{Theorem}\label{mainresult} Let $\Omega\subset\mathbb C$ be an $n$-connected domain with non-degenerate boundary components $\{\G_i\}_{i \in \{1,2...n\}}$, $z \in \Omega$.
	Then the squeezing function of $\Omega$ satisfies \begin{equation} \max_{j}\min_{i\neq j}{\arO} \leq S_{\Omega}(z)\leq \max_{j}\max_{i\neq j}{\arO}\label{formulamainresult}.\end{equation}
\end{Theorem}

\subsection{Discussion on doubly connected case} $\:$

With our notation the earlier work on the squeezing function of doubly connected domains translates to the following result. 

\begin{Theorem}\label{squeezingfunctiononannulus}
	Let $\Omega$ be a doubly connected domain with non-degenerate boundary components $\G_1,\G_2$.
	\begin{enumerate} 
		\item \cite[Theorem 1]{NTT} We have that \begin{equation}\label{formuladoublyconnected} S_{\Omega}(z)=\max\left\{r_{\Omega,1,2}(z),r_{\Omega,2,1}(z)\right\}.\end{equation} In particular, for $r\in(0,1)$ the squeezing function on the annulus $\Ar=\Delta/\overline{r\Delta}$ is given explicitly as \begin{equation} S_{\Ar}(z)=\max\left\{ |z|,\frac{\sqrt{r}}{|z|}\right\}.\label{formulaannuli}\end{equation}
		\item \cite[Theorem 2]{GR} Moreover, either one or both of the maps $\psi_{1,z}(.),\psi_{2,z}(.)$ are the unique extremal maps for this problem up to rotation.
	\end{enumerate}
\end{Theorem}

We would like to mention that \eqref{formulaannuli} follows from the explicit calculation of radii of circularly slit discs under biholomorphisms of annuli. See \cite[Lemma 3]{RW} and also \cite[Theorem 4]{NTT} for an alternative proof.

Clearly, \eqref{formuladoublyconnected} immediately follows from Theorem \ref{mainresult}. Moreover, we will remark that we can slightly modify our simple proof of Theorem \ref{mainresult} to get the second part of Theorem \ref{squeezingfunctiononannulus} as well. 

Observe that Theorem \ref{squeezingfunctiononannulus} settles the squeezing problem on doubly connected domains, as there are no conformal maps from  $\mathbb{C}^{*}:=\mathbb{C}\setminus\{0\}$ into $\Delta$ and the squeezing function of $\Omega:=\Delta^{*}=\Delta\setminus\{0\}$ is easily calculated be $S_{\Omega}(z)=|z|$, see \cite[Corollary 7.2]{DGZ1}. Therefore \eqref{formulaannuli} and the formula of the squeezing function on the punctured disc implies that the minimum of the squeezing function characterizes doubly connected domains up to biholomorphism. 

The first part of Theorem \ref{squeezingfunctiononannulus} is established in \cite{NTT} with the help of the complicated machinery of the L\"owner differential equation. It was later confirmed by \cite{GR} with potential theoretic tools and by \cite{S} with methods related to other conformal invariants. On the other hand, its uniqueness part is proven in \cite{GR} and later shown by another method in \cite{S}.  

Before we move on to the proof of Theorem \ref{mainresult} we would like to briefly discuss the higher connected case. The work of Ng, Tang and Tsai\cite{NTT} led them to conjecture that certain maps that take a domain onto a circularly slit disc are extremal mappings for the squeezing problems on higher connected domains as well. Namely with our notation they conjectured that \[S_{\Omega}(z)=\max_{j}{\min_{i\neq j}{\arO}}.\]

Their conjecture is disproven by Gumenyuk and Roth in \cite{GR} where they proved that for any $n>2$, there exists an $n$-connected domain with non-degenerate boundary $\Omega$ and $z  \in \Omega$ such that  $S_{\Omega}(z)>\max_{j}\min_{i\neq j}{\arO}$. Recently, Solynin \cite{S} showed that for any $n>2$, there exists a circularly slit disc $C$ of connectivity $n$ such that identity is the unique extremal map for the squeezing problem at $0$. 

Curiously, the circularly slit discs Solynin constructs are homogeneous in the sense that all of their slits are equidistant from the origin, whereas the examples of Gumenyuk and Roth are not. Following a question Solynin, we are curious to see if there exists a circularly slit disc with slits of different radii for which the identity mapping is an extremal map for the squeezing problem at the origin.  

We move on to the proof of our main result.

\section{Proof of Theorem \ref{mainresult}}

Throughout this section, let $\Omega$ be as in the Theorem \ref{mainresult}. We would like to start with a series of lemmas.

The following lemma simplifies the squeezing problem in one dimension. We include a proof to emphasize that functions whose images are not bounded by the unit circle cannot be extremal.

\begin{lem}\label{lemmaDGZ}
	\cite[Proposition 6.1]{DGZ1}Let $\Omega$ be an $n$-connected planar domain with non-degenerate boundary. Let $\Uot\subset\Uo$  be the set of injective holomorphic maps from $\Omega$ into $\Delta$ such that the complement of the image of $\Omega$ is relatively compact in the unit disc. Then \[ S_{\Omega}(z)=\sup_{\substack{f\in\Uot \\ f(z)=0}}\dist(0,\partial f(\Omega)).\]
\end{lem}
\begin{proof}
	Let $f:\Omega\longrightarrow\Delta$ be a conformal map such that  $f(z)=0$ and $\Delta\setminus f(\Omega)$ is not relatively compact. Note that as $\Omega$ is $n$-connected, so is $f(\Omega)$. Therefore, the condition we impose on $\Delta\setminus f(\Omega)$ holds if and only if none of the boundary components of $f(\Omega)$ is the unit circle.
	
	Let $\tilde{f(\Omega)}$ be the union of $f(\Omega)$ with compact components of its complement. Observe that $\tilde{f(\Omega)}$ is simply connected. Now let $g:\tilde{f(\Omega)}\longrightarrow\Delta$ be the Riemann map fixing $0$. Consider the function  $g \circ f \in \Uot$. We will show that $\dist(0,\partial f(\Omega))<\dist(0,\partial g(f(\Omega)))$.
	
	Let $\dc_\Omega(.,.)$ denote the inner distance obtained from the Carathéodory metric on the domain $\Omega$ and for $r\in(0,\infty)$ define $t(r)=\tanh(r/2)$. $t$ is a strictly increasing function and for all $z\in\Delta$ we have $\dist(0,z)=|z|=t(\dc_{\Delta}(0,z))$. By the conformal invariance and the monotonicity of the Carathéodory metric we have \begin{equation}\dc_{\Delta}(0,\Delta \setminus g(f(\Omega)))=
		\dc_{\tilde{f(\Omega)}}(0,\tilde{f(\Omega)}\setminus f(\Omega))>\dc_{\Delta}(0,\tilde{f(\Omega)}\setminus f(\Omega)).\label{caratheodory}\end{equation} 
	
	We note that the inequality in the right-hand side of \eqref{caratheodory} is indeed strict, since around the origin the Schwarz lemma shows that the infinitesimal version of the Carathéodory metric of $\tilde{f(\Omega)}$ is strictly greater than the infinitesimal version of the  Carathéodory metric of $\Delta$. As the function $t$ is strictly increasing we obtain \begin{equation}\dist(0,\partial g(f(\Omega)))=t(\dc_{\Delta}(0,\Delta \setminus g(f(\Omega))))>t(\dc_{\Delta}(0,f(\Omega)'\setminus f(\Omega)))=\dist(0,f(\Omega)). \end{equation}
	Hence, the lemma follows. \end{proof}

The following extremal property of conformal maps onto circularly slit discs is going to be the key lemma in our proof. 

\begin{lem}\label{lemmaRW}\cite[Theorem 3]{RW}
	Let $\Omega\subset\Delta$ be a finitely connected
	domain with outer boundary $\partial \Delta$ and let $\psi$ be the conformal mapping of $\Omega$ onto a circularly
	slit disk normalized by $\psi(0)=0$, $\psi'(0)>0$ and $\psi(\partial\Delta)=\partial\Delta$. Then $\psi'(0)\geq 1$ and $\psi'(0)=1$ if and only if $\Omega$ is a circularly slit disc, in which case $\psi={id}_{\Omega}$. 
\end{lem}

We recall the following well-known theorem of Carathéodory, for a proof see for instance \cite[Theorem 5.1.1]{KR}. 

\begin{lem}\label{lemmaC}
	Let $\Omega_1,\Omega_2$ be simply connected planar domains bounded by Jordan curves and $f:\Omega_1\longrightarrow\Omega_2$ be a biholomorphism. Then $f$ extends to a homeomorphism of $\overline{\Omega_1}$ onto $\overline{\Omega_2}$. That is there is a bijective continuous map $\tilde{f}:\overline{\Omega_1}\longrightarrow\overline{\Omega_2}$ such that $\tilde{f}|_{\Omega_1}=f$
\end{lem}

We are now ready to prove the theorem.

\textit{Proof of Theorem \ref{mainresult}.}

For $j\in\{1,2,...,n\}$ let 
\begin{equation}\label{cases}
	S_{\Omega,j}(z):=\sup_{\substack{f\in\Uo\\f(z)=0}}\{\dist(0,f(z)):f(\Gamma_j)=\partial\Delta\}.\end{equation}

Then due to Lemma \ref{lemmaDGZ} we have \begin{equation}\label{maximum} S_\Omega(z)=\displaystyle{ \max_{j\in\{1,2,...,n\}} S_{\Omega,j}(z)}. \end{equation}

We claim that we have
\begin{equation}\displaystyle{\min_{i\neq j}} \arO \leq S_{\Omega,j}(z)\leq \displaystyle{\max_{i\neq j}} \arO.\label{claimed}\end{equation}

The left-hand side of \eqref{claimed} follows trivially by observing that biholomorphisms onto circularly slit discs are candidates for the supremum given in \eqref{cases}.

For the right-hand side, note that by a standard normal families argument we can easily see that the supremum in \eqref{cases} is always achieved. Thus, without loss of generality, we assume that $\Omega$ is a domain such that the identity function is an extremal mapping for \eqref{cases}. We'll show that $$\dist(0,\partial \Omega)\leq\max_{i \neq j}\arO.$$

Let $\psi:=\psi_{z,j}$ be the conformal map of $\Omega$ onto a circularly slit disc $C:=C_{z,j}(\Omega)$ fixing $0$ and $\partial\Delta$. Without loss of generality assume that ${\psi}'(z)>0$. 

Consider the function $g(z):={z}/{{\psi}^{-1}(z)}$ on $C$. Here, ${\psi}^{-1}$ denotes the inverse of $\psi$. As ${\psi}^{-1}(z)=(\psi^{-1})'(0)z+\mathcal{O}(z^2)$ we have 
\[\displaystyle{\lim_{z \to 0} g(z)}=\displaystyle{\lim_{z \to 0} \dfrac{z}{{\psi}^{-1}(z)}=\dfrac{1}{({\psi}^{-1})'(0)}}={\psi}'(0).\]

We conclude that $g$ is bounded near $0$ so it extends to $C$ analytically. With an abuse of notation, we use $g$ to denote its analytic extension. 

Due to Lemma \ref{lemmaRW} ${\psi}'(0)\geq1$, therefore we must have that $g(0)\geq1$. If $g(0)=1$ then the mapping $\psi$ is identity. This shows that $\Omega$ is a circularly slit disc so \eqref{claimed} follows. 

Now suppose $g(0)>1$. In this case, $\Omega$ is not a circularly slit disc. 

For $k\in\mathbb{N}$ let $\varepsilon>0$ be small enough so that the curves $\gamma_{i,k}:=\{z\in C:\dist(0,f(\G_i))=\frac{\varepsilon}{k+1}\}$ are disjoint Jordan curves in $C$. As ${\psi}^{-1}$ is injective, their preimages are Jordan curves as well. Let $C_k$ be the domain bounded by $\partial\Delta$ and $\gamma_{i,k}$ for $k\in\mathbb{N}$ and $i\neq j\in\{1,2...n\}$.

Let $p\in\partial C_k$ be arbitrary and choose a small enough neighbourhood around $p$, $U_p$ so that $U_p\cap C_k$ is bounded by a Jordan curve so that the image of $U_p\cap C_k$ under $\psi^{-1}$ is bounded by a Jordan curve as well. It therefore follows from Lemma \ref{lemmaC} that $\psi^{-1}$ extends continuously to $p$. As $p$ was arbitrary we conclude that $\psi^{-1}$ extends continuously to $\overline{C_k}$, hence so does $g$.  

The maximum principle asserts that \begin{equation}\label{principle}
	\sup_{x\in \overline{C_k}}|g(x)|=\sup_{x\in\partial C_k}|g(x)|.\end{equation}

As $\Omega$ is not a circularly slit disc, $g(0)>1$. Thus, by \eqref{principle} we must have a sequence $z_k \in \partial C_k$ with $|g(z_k)|>1$. Noting that $\displaystyle{\lim_{x\longrightarrow z}|g(x)|} = 1$ for all $z\in\partial\Delta$ we can see that $\{z_k\}\in\partial C_k \setminus \partial \Delta$. Moreover \begin{equation}1 < |g(z_k)|  \leq \dfrac{\max_{\partial C_k \setminus \partial \Delta}|z|}{|{\psi}^{-1}(z_k)|}. \label{nearlyfinished}\end{equation}

From \eqref{nearlyfinished} we can directly see that \begin{equation}|{\psi}^{-1}(z_k)|<\max_{\partial C_k \setminus \partial \Delta}|z|=\dfrac{\varepsilon}{k+1}+\max_{i\neq j}\arO. \label{nearlyfinished2}\end{equation}

It is clear from the definition of $C_k$ that any convergent subsequence of $\{z_k\}$ must converge to $\partial C$. From \eqref{nearlyfinished2} we obtain
\begin{equation}\displaystyle{\limsup_{k\rightarrow\infty} |{\psi}^{-1}(z_k)|}\leq \displaystyle{\limsup_{k \to \infty} \left(\dfrac{\varepsilon}{k+1}+\max_{i\neq j}\arO\right)}= \max_{i\neq j}\arO. \label{nearlyfinished3}\end{equation}

\eqref{nearlyfinished3} shows that $\dist(0,\partial\Omega)\leq\max_{i\neq j}\arO$. So we have \eqref{claimed}. 

By \eqref{maximum}, the theorem follows. $\endofproof$

Before we conclude, we give our final observations. 

If all of the slits of the circularly slit disc $C$ given above are not on the same circle, then we see that the inequality on the right-hand side of \eqref{claimed} is strict. To see this suppose that $\Omega$ is a domain such that the identity function is the extremal map for \eqref{cases} for $0\in\Omega$. If $\Omega$ is a circularly slit disc, then clearly the inequality is strict. If not then at $0$, the norm of the derivative of the map $\psi$ given above is strictly more than $1$ due to Lemma \ref{lemmaRW}. As in the proof of Theorem \ref{mainresult}, suppose that $$\left| \frac{1}{{({\psi}^{-1})}'(0)} \right|=|\psi'(0)|=\alpha>1$$ Defining $g(z)=\dfrac{2z}{(\alpha+1){\psi}^{-1}(z)}$ and applying the same argument shows that $$\dist(0,\partial\Omega)\leq\dfrac{2}{(\alpha+1)}\max_{i\neq j}\arO<\max_{j\neq i}\arO$$ so we have the desired conclusion.

On the other hand, if $\Omega$ is a circularly slit disc with slits lying on the same circle, the argument we gave above shows that identity is the unique extremal domain for squeezing problem for $0\in\Omega$ fixing $\partial \Delta$. Therefore we recover \cite[Lemma 1]{S}.

Moreover in the doubly connected case, circularly slit discs have only one slit, so all of their slits automatically lie on the same circle. Following the arguments above, we can also get the second part of Theorem \ref{squeezingfunctiononannulus}. 

$\:$

$\:$

\textbf{Acknowledgements.} This note was written during the author's studies at Middle East Technical University. The author is grateful to the mathematics department and his M.Sc. advisor {\"O}zcan Yaz{\i}c{\i}.

{}

\end{document}